\documentclass{au}
\usepackage{amssymb,amsfonts,amscd}
\theoremstyle{plain}
\newtheorem{thm}{Theorem}[section]

\newtheorem{prop}[thm]{Proposition}

\theoremstyle{definition}
\newtheorem*{ex}{Example}

\newtheorem{defi}[thm]{Definition}

\newtheorem*{rem}{Remark}

\begin{document}

\author{E. Graczy\'{n}ska}
\address{Technical University of Opole, Institute of Mathematics
\newline
ul. Luboszycka 3, 45-036 Opole, POLAND} \email{egracz@po.opole.pl
\hspace{1cm} http://www.egracz.po.opole.pl/}
\author{D. Schweigert}
\address{Technische Universit\"{a}t Kaiserslautern, Fachbereich Mathematik
\newline
Postfach 3049 \\ 67653 Kaiserslautern, Germany}
\email{schweige@mathematik.uni-kl.de}

\title{M-hyperquasivarieties}

\begin{abstract}
We consider the notion of M-hyper-quasi-identities and
M-hyperquasivarieties, as a common generalization of the concept of
quasi-identity (hyper-quasi-identity) and  quasivariety
(hyper-quasivariety) invented by A. I. Mal'cev, cf. \cite{AIM}, cf.
\cite{VAG} and hypervariety invented by the authors in \cite{9},
\cite{4} and hyperquasivariety \cite{EGDS3}. The results of this
paper were presented on the 69th Workshop on General Algebra, held
at Potsdam University (Germany) on March 18-20, 2005.


\end{abstract}

\maketitle

\section{Notations}

An identity is a pair of terms where the variables are bound by
universal quantifiers. Let us take the following medial identity
as an example
\[
\forall u\forall x\forall y\forall w\,(u\cdot x)\cdot (y\cdot
w)=(u\cdot y)\cdot (x\cdot w).
\]
Let us look at the following hyperidentity
\[
\,\forall F\forall u\forall x\forall y\forall
w\,F(F(u,v),F(x,y))=F(F(u,x),F(v,y)).
\]
The hypervariable $F$ is considered in a very specific way.
Firstly every hypervariable is restricted to functions of a given
arity. Secondly $F$ is restricted to term functions of the given type.
Let us take the variety $A_{n,0}$ of abelian groups of finite exponent
$n.$ Every binary term $t\equiv t(x,y)$ can be presented by
$t(x,y)=ax+by$ with $a,b\in \mathbb N_{0}.$ If we substitute the
binary hypervariable $F$ in the above hyperidentity by $ax+by,$
leaving its variables unchanged, we get
\[
a(au+bv)+b(ax+by)=a(au+bx)+b(av+by).
\]
This identity holds for every term $t(x,y)=ax+by$ for the variety
$A_{n,0}$. Therefore we say that the hyperidentity holds for the
variety $A_{n,0}.$

{\bf M-hyper-quasi-identities}

In the sequel we use the definition of {\em hyperterm} from
\cite{9} and \cite{4}. We accept the notation of \cite{PMC}, \cite{5},
\cite{VAG}, \cite{DHRK} and \cite{JP}.

Given a monoid M of hypersubstitutions of a given type $\tau$.

We recall only the definitions of \cite{JP} of the fact that a hyperidentity
is satisfied in an algebra as an M-hyperidentity of a given type and the
notion of M-hypervariety, \cite{EG}:

\begin{defi}\label{D:2.1}
An algebra ${\bf A}$ satisfies a hyperidentity $h_{1}=h_{2}$ as
M-hyperidentity if for every M-hypersubstitution substitution of the
hypervariables by terms (of the same arity) of ${\bf A}$ leaving the
variables unchanged, the identities which arise hold in ${\bf A}$.
In this case, we write ${\bf A} \models_{M}^{H} (h_{1}=h_{2})$. A
variety $V$ satisfies a hyperidentity $h_{1}=h_{2}$ as
M-hyperidentity if every algebra in the variety does; in this case,
we write $V\models_{M}^{H} (h_{1}=h_{2})$.
\end{defi}
\begin{defi}\label{D:2.2}
A class $V$ of a algebras of a given type is called an M-hypervariety if
and only if $V$ is defined by a set of M-hyperidentities.
\end{defi}

The following was proved in \cite{EG}:

\begin{thm} \label{T:2.3}
A variety $V$ of type $\tau$ is defined by a set $\Sigma$ of
M-hyperidentities if and only if $V = HSPD_{M}(V)$, i.e. $V$ is a
variety closed under M-derived algebras of type $\tau$. Moreover, in
this case, the set $\Sigma$ is then M-hypersatisfied in $V$ and $V$
is the class of all M-hypermodels of $\Sigma$, i.e. $V
=MHMod(\Sigma)$.
\end{thm}
We recall from \cite{AIM} and \cite{VAG}:
\begin{defi}\label{D:2.4}
A quasi-identity $e$ is an implication of the form:
\begin{center}
$(t_{0} = s_{0}) \wedge ... \wedge (t_{n-1} = s_{n-1}) \rightarrow (t_{n} = s_{n})$.
\end{center}
where $t_{i} = s_{i}$ are $k$-ary identities of a given type, for $i = 0,...,n$.

A qausi-identity above is {\em satisfied in an algebra} ${\bf A}$ of a given type
if and only if the following implication is satisfied in ${\bf A}$:  given a
sequence $a_{1},...,a_{k}$ of elements of $A$. If this elements
satisfy the equations $t_{i}(a_{1},...,a_{k}) = s_{i}(a_{1},...,a_{n})$
in ${\bf A}$, for $i = 0, 1,..., n-1$, then the equality
$t_{n}(a_{1},...,a_{k}) = s_{n}(a_{1},...,a_{k})$ is satisfied in ${\bf A}$.
In that case we write:
\begin{center}
${\bf A} \models (t_{0} = s_{0}) \wedge ... \wedge (t_{n-1} = s_{n-1}) \rightarrow (t_{n} = s_{n})$.
\end{center}

A quasi-identity $e$ is {\em satisfied in a clas}s $V$ of algebras of a given type,
if and only if it is satisfied in all algebras ${\bf A}$ belonging to $V$.
In that case we write:
\begin{center}
$V \models (t_{0} = s_{0}) \wedge ... \wedge (t_{n-1} = s_{n-1}) \rightarrow (t_{n} = s_{n})$.
\end{center}

\end{defi}

Following the ideas of \cite{9}, part 5, cf. \cite{2}, we modify the
definition above in the following way:
\begin{defi}\label{D:2.5}
A hyper-quasi-identity $e^{*}$ is an implication of the form:
\begin{center}
$(T_{0} = S_{0}) \wedge ... \wedge (T_{n-1} = S_{n-1}) \rightarrow (T_{n} = S_{n})$.
\end{center}
where $T_{i} = S_{i}$ are hyperidentities of a given type, for $i = 0,...,n$.

A hyper-quasi-identity $e^{*}$  is M-hyper-satisfied (holds) in an algebra
${\bf A}$ if and only if the following implication is satisfied:\\
If $\sigma$ is a hypersubstitution of $M$ and the elements
$a_{1},...,a_{n} \in A$ satisfy the equalities
$\sigma(T_{i})(a_{1},...,a_{k}) = \sigma(S_{i})(a_{1},...,a_{k})$ in
${\bf A}$, for $n = 0,1,...,n-1$, then the equality
$\sigma(T_{n})(a_{1},...,a-{k}) = \sigma(S_{n})(a_{1},...,a_{k})$
holds in ${\bf A}$.

We say then, that $e^{*}$ is an M-hyper-quasi-identity of ${\bf A}$ and write:
\begin{center}
${\bf A} \models^{H}_{M} (t_{0} = s_{0}) \wedge ... \wedge (t_{n-1} = s_{n-1}) \rightarrow (t_{n} = s_{n})$.
\end{center}

A hyper-quasi-identity $e^{*}$  is M-hyper-satisfied (holds) in a
class $V$ if and only if it is M-hypersatisfied in any algebra of
$V$.

\end{defi}
By other words, M-hyper-quasi-identity is a universally closed  Horn
$\forall x \forall \sigma$-formulas, where x vary over all sequences of
individual variables (occuring in terms of the implication) and $\sigma$
vary over all hypersubstitutions of $M$. Our modification coincides
with Definition 5.1.3 of \cite{9} (cf. Definition 2.3 of \cite{CCKD}).

\begin{rem}
All hyper-quasi-identities and hyperidentities are written without
quantifiers but they are considered as universally closed Horn
$\forall$-formulas (cf. \cite{AIM}). In case of a trivial monoid
$M$, the notion of M-hypersatisfaction reduces to the notion of
classical satisfaction of \cite{GB}, \cite{PMC}. If $M$ is the
monoid of all hypersubstiutions of a given type $\tau$, then the
notion of M-hyperidentity and M-hyperquasiidentity reduces to the
hyperidentity (cf. \cite{4}) and hyperquasi-identity(cf.
\cite{EGDS3}).
\end{rem}

\section{Examples of hyper-quasi-identities}
\subsection{Hyper-quasi-identities for abelian algebras}
\begin{defi}\label{D:3.1}
An algebra ${\bf A}$ is called abelian if for every $n > 1$ and every n-ary term
operation $f$ of ${\bf A}$ and for all $u,v,x_{1},...,x_{n-1},y_{1},...,y_{n-1}$
the following equivalence holds:
\begin{center}
$f(u,x_{1},...,x_{n}) = f(u,y_{1},...,y_{n}) \leftrightarrow f(v,x_{1},...,x_{n}) = f(v,x_{1},...,x_{n})$
\end{center}
A variety $V$ is called abelian, if each algebra of $V$ is abelian.
\end{defi}
It follows from  \cite[p. 40]{DHRK}, \cite[p. 290]{2}:

\begin{prop}\label{P:3.1}
An algebra ${\bf A}$ is abelian if and only if the following  hyper-quasi-identity
holds in ${\bf A}$:
\begin{center}
$F(u,x_{1},...,x_{n}) = F(u,y_{1},...,y_{n}) \leftrightarrow
F(v,x_{1},...,x_{n}) = F(v,y_{1},...,y_{n})$.
\end{center}
\end{prop}
\begin{ex}
The variety $RB$ of rectangular bands fulfills the hyperidentities of type
$(1,2,3,...,n,...)$:\\
$F(x,...,x) = x$, $F(F(x_{11},...,x_{1n}),x_{2},...,x_{n}) = F(x_{11},x_{2},...,x_{n})$, \\
$F(x_{1},...,x_{n-1},F(x_{n1},...,x_{nn})) = F(x_{1},...,x_{n-1},x_{nn})$.\\
We can derive by the associative hyperidentity:
$F(x,F(y,z)) = F(F(x,y),z)$, that the following hyper-quasi-identity holds in $RB$:\\
$F(u,x_{1},...,x_{n}) = F(u,y_{1},...,y_{n}) \rightarrow
F(v,u,x_{1},...,x_{n-1}) = F(v,u,y_{1},...,y_{n-1}) \rightarrow
F(v,x_{1},...,x_{n}) = F(v,y_{1},...,y_{n})$, i.e. the variety  $RB$ is abelian.
\begin{rem}
For further examples of abelian varieties consult \cite{DHRK} and \cite{11}.
\end{rem}

\end{ex}
\subsection{Semidistributive lattices}
A lattice is {\em joinsemidistributive} if it satisfies the following
condition, cf. \cite[p. 82]{DHRK},  \cite[p. 141]{VAG}:
\begin{center}
($SD_{\vee}$) $x \vee y = x \vee z$ implies $x \vee y = x \vee (y \wedge z)$
\end{center}
The {\em meetsemidistributivity} is defined by duality:
\begin{center}
($SD_{\wedge}$) $x \wedge y = x \wedge z$ implies $x \wedge y = x \wedge (y \vee z)$
\end{center}

A lattice is {\em semidistributive} if it is simultaneously join and meet
semidistributive.

\begin{prop}\label{P:3.2}
Let ${\bf L} = (L, \wedge, \vee)$ be a lattice which is semidistributive.
Then the following hyper quasi-identity is hypersatisfied in ${\bf L}$:
\begin{center}
$(F(x,y) = F(x,z)) \rightarrow (F(x,y) = F(x,G(y,z)))$
\end{center}
\end{prop}
\begin{proof}
We consider all cases on hypervariables of a semidistributive lattice.

Case 1. $F(x,y) := x$.\\
Obviously $(x=x) \rightarrow (x =x)$.

Case 2. $F(x,y) := y$. Consider $G(y,z) := y, z, y \wedge z, y \vee z$.\\
Then the following quasi-identities are satisfied in ${\bf L}$: \\
$(y = z) \rightarrow (y = y)$, $(y =z) \rightarrow (y = z)$,
$(y = z) \rightarrow (y = y \wedge z)$,\\
$(y = z) \rightarrow (y = y \vee z)$.

Case 3. $F(x,y) := x \wedge y$. Consider $G(y,z) := y,z, y\wedge z, y \vee z$. \\
Then the following quasi-identities hold in ${\bf L}$: \\
$(x \wedge y = x \wedge z) \rightarrow  (x \wedge y = x \wedge y)$,
$(x \wedge y = x \wedge z) \rightarrow (x \wedge y = x \wedge z)$, \\
$(x \wedge y = x \wedge z) \rightarrow (x \wedge y = x \wedge (y\wedge z))$,
$(x \wedge y = x \wedge z) \rightarrow (x \wedge y = x \wedge (y \vee z))$.
The last one follows from the  meet semidistributivity of ${\bf L}$.

Case 4.  $F(x,y) := x \vee y$. Consider $G(y,z) := y,z, y\wedge z, y \vee z$.\\
Then the following quasi-identities hold in ${\bf L}$:\\
$(x \vee y = x \vee z) \rightarrow  (x \vee y = x \vee y)$,
$(x \vee y = x \vee z) \rightarrow (x \vee y = x \vee z)$, \\
$(x \vee y = x \vee z) \rightarrow (x \vee y = x \vee (y\wedge z))$,
$(x \vee y = x \vee z) \rightarrow (x \vee y = x \vee (y \vee z))$.
The prelast one follows from the  join semidistributivity of ${\bf L}$.
\end{proof}

\section{M-derived algebras}
We recall from \cite{9}, \cite{4}, \cite{EG}, cf. \cite[p. 145]{PMC}
the notion of a {\em M-derived algebra} and the {\em M-derived
class} of algebras. Given a type $\tau =
(n_{0},n_{1},...,n_{\gamma},...)$. An algebra ${\bf B}$ of type
$\tau$ is called an M-derived algebra of ${\bf A} =
(A,f_{0},f_{1},...,f_{\gamma},...)$ if there exist term operations
$t_{0},t_{1},...,t_{\gamma},...$ of type $\tau$ such that ${\bf B} =
{\bf A}^{\sigma} = (A,t_{0},t_{1},...,t_{\gamma},...)$ with
$t_{\gamma}$ being equal to the term
$\sigma(f_{\gamma}(x_{0},...,x_{\tau(f)-1}))$ for a fixed $\sigma
\in M$.

For a class $K$ of algebras of type $\tau$ we denote by ${\bf D_{M}}(K)$ the
class of all M-derived algebras (of type $\tau$) of $K$. \\
A class $K$ is called M-{\em deriverably closed} if and only if
${\bf D_{M}}(K) \subseteq K$, for a given monoid $M$ of hypersubstitutions of
type $\tau$.

For a given algebra ${\bf A}$ we denote by $QId({\bf A})$ and
$MHQId({\bf A})$ the set of all quasi-identities and
M-hyper-quasi-identities satisfied (M-hypersatisfied) in ${\bf A}$,
respectively.  Similarly for a class $K$, $QId(K)$ and $MHQId(K)$
denote the set of all quasi-identities and M-hyper-quasi-identities
satisfied (M-hypersatisfied) in $K$, respectively. Following
\cite{VAG} by ${\bf Q}(K)$ and $MH{\bf Q}(K)$ we denote the class of
all algebras of a given type satisfying (M-hypersatisfying) all the
quasi-identities and M-hyper-quasi-identities of $K$, respectively.
The transformation $T$ of an quasi-identity $e$ into
hyper-quasi-identity $e^{*}$ is defined in a natural way. Similarly
$T^{-1}$ transforms every hyper-quasi-identity $e^{*}$ into the
quasi-identity $e$. (cf. \cite{4}).

\begin{prop}\label{P:4.1}
Given a class $K$ of algebras of type $\tau$. Then the following equality
holds:
\begin{center}
$MHQId(K) = T(QId({\bf D_{M}}(K)))$.
\end{center}
\end{prop}
\begin{proof}
To prove the inclusion $\subseteq$, given a hyper-quasi-identity
$e^{*}$ of $K$ and an algebra an M-derived algebra ${\bf B} = {\bf
A}^{\sigma}$ for ${\bf A} \in K$ and some $\sigma \in M$. Then by
definition 1.5 $T^{-1}(e)$ is satisfied in ${\bf B}$ as a
quasi-identity, i.e. $e^{*} \in T(QId({\bf D_{M}}(K)))$.

For a proof of the converse inclusion, given a quasi-identity $e$
satisfied in ${\bf D_{M}}(K)$ and an algebra ${\bf A}$ of $K$. Let
$a_{1},...,a_{k} \in A$. Given an hypersubstitution $\sigma$ of  $M$
and consider $\sigma(e)$ in ${\bf A}$. As ${\bf A} \in {\bf
D_{M}}(K)$, then $\sigma(e)$ may be considered as a quasi-identity
of the M-derived algebra ${\bf A}^{\sigma}$ (cf. \cite{10}).\\
Assume that $p_{i}^{{\bf A}^{\sigma}}(a_{1},,,a_{k}) = q_{i}^{{\bf
A}^{\sigma}}(a_{1},...,a_{k})$, therefore as ${\bf A} \in {\bf
D_{M}}(K)$ we obtain that $p_{n}^{{\bf A}^{\sigma}}(a_{1},...,a_{k})
= q_{n}^{{\bf A}^{\sigma}}(a_{1},...,a_{n})$, i.e. that the $T(e)$
holds in $K$ as an M-hyper-quasi-identity, i.e. $T(e) \in MHQId(K)$.
\end{proof}
The role of derived algebras in solvability questions may be visualized by the
following:
\begin{thm}\label{T:4.2}
Given a locally finite hypervariety $V$. The class $K$ of all locally solvable
algebras in $V$ is a hypervariety.
\end{thm}
\begin{proof}
By corollary 7.6, of \cite{DHRK} the class of all locally solvable algebras in a locally
finite variety $V$ is a variety. Due to our theorem 2.3 we need to show that
$K$ is deriverably closed. Given a derived algebra ${\bf B}$ of an algebra
${\bf A} \in K$ and its finite subalgebra ${\bf F}^{*}$. Then ${\bf F}^{*}$
is a derived subalgebra of a finite subalgebra ${\bf F}$ of ${\bf A}$, i.e.
${\bf F}^{*} = {\bf F}^{\sigma}$, for a hypersubstitution $\sigma$.
From theorem 5.7 of \cite{10} we conclude, that as ${\bf F}$ is solvable,
therefore ${\bf F}^{\sigma}$ is solvable.
\end{proof}

\section{M-hyper-quasi-varieties}
We reformulate the notion of {\em quasivariety} invented by A. I.
Mal'cev in \cite[p. 210]{AIM} and {\em hyperquasivariety} of
\cite{EGDS3} for the case of {\em M-hyper-quasivariety} of a given
type in a natural way:
\begin{defi}\label{D:5.1}
A class $K$ of algebras of type $\tau$ is called an M-hyperquasivariety if there
is a set  $\Sigma$ of M-hyper-quasi-identities of type $\tau$ such that
$K$ conisits exactly of those algebras of type $\tau$ that M-hypersatisfy all
the hyper-quasi-identities of $\Sigma$.
\end{defi}
From proposition 3.1 we obtain immediately:
\begin{thm}\label{T:5.2}
A quasivariety $K$ of algebras given type is a M-hyperquasivariety if and
only if it is  M-deriverably closed.
\end{thm}
\vspace{1cm}
In the sequel we use the standard notation: ${\bf S}$ for the operator of
creating \emph{subalgebras}, ${\bf P}$, ${\bf P_{s}}$, ${\bf P_{r}}$, ${\bf P_{u}}$
and ${\bf P{\omega}}$ for \emph{products}, \emph{subdirect products},
\emph{reduced products}, \emph{ultraproducts} and \emph{direct products of finite families of
structures}, respectively.  ${\bf L}$ and ${\bf L_{s}}$ will denote the operators of
\emph{direct limits} and \emph{superdirect limits}.

Following  A. I. Mal'cev \cite[p. 153, 215]{AIM} for a given class $K$
of algebras of a given type, by ${\bf S}(K)$, ${\bf P}(K)$, ${\bf P_{s}}(K)$, ${\bf P_{r}}(K)$,
${\bf P_{u}}(K)$ we denote the class of algebras isomorphic to all possible subalgebras, direct products,
subdirect products, reduced (filtered) products or ultraproducts of algebras
of $K$, respectively. ${\bf P_{\omega}}$ is the class of algebras
isomorphic to the direct products of finite families of structures of $K$.
Similarly, ${\bf L}(K)$ and ${\bf L_{s}}(K)$ will be the class of algebras isomorphic to
direct and superdirectlimits of algebras of $K$, respectively (cf. \cite[p.21]{VAG}.

Adding a trivial system to $K$ we obtain the class $K_{0}$.

We reformulate the result of \cite{4}:
\begin{prop}\label{P:5.3}
Given a class $K$ of algebras and a monoid $M$ of hypersubstitutions,
then the following inclusions holds:\\
1) ${\bf D_{M}}{\bf S}(K) \subseteq {\bf S}{\bf D_{M}}(K)$; \\
2) ${\bf D_{M}}{\bf P}(K) \subseteq {\bf P}{\bf D_{M}}(K)$;\\
3) ${\bf D_{M}}{\bf P_{\omega}}(K) \subseteq {\bf P_{\omega}}{\bf D_{M}}(K)$;\\
4) ${\bf D_{M}}{\bf P_{s}}(K) \subseteq {\bf P_{s}}{\bf D_{M}}(K)$;\\
5) ${\bf D_{M}}{\bf P_{r}}(K) \subseteq {\bf P_{r}}{\bf D_{M}}(K)$;\\
6) ${\bf D_{M}}{\bf P_{u}}(K) \subseteq {\bf P_{u}}{\bf D_{M}}(K)$;\\
7) ${\bf D_{M}}{\bf L}(K) \subseteq {\bf L}{\bf D_{M}}(K)$;\\
8) ${\bf D_{M}}{\bf L_{s}}(K) \subseteq {\bf L_{s}}{\bf D_{M}}(K)$.
\end{prop}
\begin{proof}
Obviously an isomorphism respects the inclusions above.  \\
The inclusions 1) and 2) were proved in \cite{4}.\\
3) and 4) immediately follows from 2). \\
To show 5) assume that an algebra ${\bf B} \in {\bf D_{M}}{\bf P_{r}}(K)$.\\
Let ${\bf B}= {\bf A}^{\sigma} = (A, t^{\bf A}_{0},t^{\bf A}_{1},...,t^{\bf A}_{\gamma},...)$
is an M-derived algebra of a reduced product ${\bf A} = (A, f^{\bf A}_{0},f^{\bf A}_{1},...,f^{\bf A}_{\gamma},...)$,
where $A = (\Pi_{i \in I}A_{i})/\sim_{F}$, for a set $I$, a family
$({\bf A}_{i})_{i \in I}$ of algebras of type $\tau$ from $K$ and a filter $F$ over $I$,
i.e. ${\bf A} = (\Pi_{i \in I}{\bf A}_{i})/F$, where for any $n$-ary functional
symbol $f$ the following holds  (cf. \cite[p. 13]{VAG}):
\begin{center}
$f^{\bf A}(a_{0}/F,...,a_{n-1}/F) = a_{n}/F$ if and only if
$\{i \in I: f^{{\bf A}_{i}}(a_{0},...,a_{n-1}) = a_{n} \} \in F$.
\end{center}
Note that by the induction on the complexity of a term $t$ of type $\tau$
one may show, that the following holds for any $n$-ary polynomial $t$ of type
$\tau$:
\begin{center}
$t^{\bf A}(a_{0}/F,...,a_{n-1}/F) = a_{n}/F$ if and only if
$\{i \in I: t^{{\bf A}_{i}}(a_{0},...,a_{n-1}) = a_{n} \} \in F$.
\end{center}
From the above equality it follows that the algebra ${\bf B} =
\Pi_{i \in I}{\bf A}_{i}^{\sigma}/F$, where ${\bf A}_{i}^{\sigma} =
(A_{i}, t_{0}^{{\bf A}_{i}},t_{1}^{{\bf A}_{i}},...,t_{\gamma}^{{\bf
A}_{i}},...)$, for $i \in I$ is an M-derived algebra of ${\bf
A}_{i}$, i.e. ${\bf B} \in {\bf P_{r}}{\bf D_{M}}(K)$.

The inclusion 6) is an immediate consequence of 5).

To prove 7), let ${\bf B} \in {\bf D_{M}}{\bf L}(K)$, i.e. given an M-derived algebra
${\bf B}$ of an algebra ${\bf A}$, where ${\bf A} = (A,f_{0}^{\bf A},f_{1}^{\bf A},...,f_{\gamma}^{\bf A},...)$ is a direct limit of
a triple $\Lambda = (I, {\bf A}_{i},g_{ij})$, where $I = (I, \leq)$ is an
up-directed set, $({\bf A}_{i})_{i \in I}$ is a family of algebras of a
given type $\tau$ and $\{g_{ij} : i,j \in I , i \leq j\}$ is a family of
homomorphisms of ${\bf A}_{i}$ into ${\bf A}_{j}$ called a \emph{direct
spectrum} over $(I, \leq)$, cf. \cite[p. 17]{VAG}. For a given direct spectrum
$\Lambda = ({\bf A}_{i},g_{ij})$, we consider the quotient set
$A = \bigcup _{i \in I} A_{i} \times \{i\}/\equiv$, where
\begin{center}
$(a,i) \equiv (b,j)$ if and only if $(\exists k \in I)(i,j \leq k, g_{ik}(a) = g_{jk}(b))$.
\end{center}
Let $<a,i>$ denotes the equivalence class by $\equiv$ containing $(a,i)$.
The operations on $A$ are defined by setting for any operation symbol of type $\tau$:
\begin{center}
$f^{\bf A}(<a_{0},i_{0}>,...,<a_{n-1},i_{n-1}>) =
<f^{{\bf A}_{j}}(g_{i_{0}j}(a_{0}),...,g_{i_{n-1}j}(a_{n-1})),j>$.
\end{center}
Note, that as $g_{ij}$ are homomorphisms for all $i \leq j$, $i,j \in I$,
therefore for any polynomial symbol $p$ of type $\tau$, the following holds:
\begin{center}
$p^{\bf A}(<a_{0},i_{0}>,...,<a_{n-1},i_{n-1}>) =
<p^{{\bf A}_{j}}(g_{i_{0}j}(a_{0}),...,g_{i_{n-1}j}(a_{n-1})),j>$.
\end{center}
As ${\bf B} = {\bf A}^{\sigma}$ for a hypersubstitution $\sigma$ of $M$,
therefore ${\bf B} \in {\bf L}{\bf D_{M}}(K)$, namely ${\bf B}$ is a direct limit of the
triple $\Lambda^{\sigma} = (I,{\bf A}_{i}^{\sigma},g_{ij})$ of M-derived algebras of
$K$, i.e. of a derived spectrum of $\Lambda$.

To prove 8), recall \cite[17]{VAG}, that a direct spectrum
$\Lambda = (I,{\bf A}_{i},gij)$ is called \emph{superdirect} if the mappings
$g_{ij}: A_{i} \rightarrow A_{j}$ are surjective. The family
$({\bf A}_{i})_{i \in I}$ is referred to as \emph{superdirect family}.
The direct  limit of a superdirect spectrum is called the \emph{superdirect limit}.
Note, that the derived spectrum $\Lambda^{\sigma}$ of  a superdirect spectrum
$\Lambda$ is superdirect. This together with 7) proves 8).
\end{proof}
Via corollary 2.3.4 \cite[p. 79]{VAG} and proposition 5.1 we obtain another
characterization of hyperquasivarieties:
\begin{prop}\label{P:5.4}
For any class $K$ of algebras, the following assertions hold:\\
(1) $MH{\bf Q}(K) = {\bf S}{\bf P}{\bf P_{u}}{\bf D_{M}}(K)$;\\
(2) $MH{\bf Q}(K) = {\bf S}{\bf P_{u}}{\bf P}{\bf D_{M}}(K)$;\\
(3) $MH{\bf Q}(K) = {\bf S}{\bf P_{u}}{\bf P_{\omega}}{\bf D_{M}}(K)$;\\
(4) $MH{\bf Q}(K) = {\bf S}{\bf L_{s}}{\bf P}{\bf D_{M}}(K) = {\bf L_{s}}{\bf S}{\bf P}{\bf D_{M}}(K) = {\bf L_{s}}{\bf P_{s}}{\bf S}{\bf D_{M}}(K)$.
\end{prop}

Via proposition 4.1 we obtain the following reformulation of Mal'cev
theorems of \cite[p. 214 - 215]{AIM}:
\begin{thm}\label{T:5.5}
A class $K$ of algebras of a given type is a hyperquasivariety if and only if
$K$\\
i) is ultraclosed;\\
ii) is heraditery;\\
iii) is multiplicatively closed; \\
iv) contains a trivial system; \\
v) is M-deriverably closed.
\end{thm}
\begin{thm}\label{T:5.6}
For every class $K$ of algebras of a given type we have
\begin{center}
$MH{\bf Q}(K) = {\bf S}{\bf P_{r}}{\bf D_{M}}(K_{0})$
\end{center}
\end{thm}
\begin{rem}
If we accept, that the direct product of an empty family of algebras of a given
type  is a trivial algebra of a given type, we may remove condition iv) of
theorem  5.4 and substitute $K_{0}$ by $K$ in theorem 5.4.
(cf. Corollary 2.3. of \cite[p. 78]{VAG}).
\end{rem}

\end{document}